\pdfoutput=1
\documentclass[11pt,oneside]{amsart}
\usepackage[leqno]{amsmath}
\usepackage{amssymb,amsthm} 
\usepackage[parfill]{parskip} 
\usepackage{libertine} 
\usepackage[T1]{fontenc}
\usepackage{textcomp}
\usepackage[varqu,varl]{inconsolata} 
\usepackage{amsmath,amsthm}
\usepackage[libertine,bigdelims,vvarbb]{newtxmath}
\usepackage[supstfm=libertinesups,%
  supscaled=1.2,%
  raised=-.13em]{superiors}
\usepackage{graphicx} 
\usepackage{calc}

\newcommand{\beq}{\begin{equation}}
\newcommand{\eeq}{\end{equation}}
\numberwithin{equation}{section} 
\newtheorem{thm}[equation]{Theorem}
\newtheorem{lem}[equation]{Lemma}

\theoremstyle{remark}

\theoremstyle{definition}
\newtheorem{defn}[equation]{Definition}

\title{On the difficulty of proving P equals NP in ZFC}
\usepackage[paperwidth=7in, paperheight=10in, textwidth=4.75in, textheight=8in, includehead=true]{geometry}
\usepackage{url}
\usepackage[pdfborderstyle={/S/U/W 1}]{hyperref}

\def\setmax{\operatorname{setmax}}

\let\emptyset\emptysetAlt
\hypersetup{
pdftitle={ZFC},
pdfauthor={\textcopyright\ S. G. Williamson}
}
\author{S. Gill Williamson}
\thanks{Department of Computer Science and Engineering, 
University of California San Diego; \url{http://cseweb.ucsd.edu/~gill/}; 
gwilliamson@ucsd.edu;
}

\begin{document}
\begin{abstract}
In Friedman's 1998 Annals of Mathematics paper~\cite{hf:nlc} and his downloadable manuscript~\cite{hf:alc},  he presents numerous combinatorial statements and shows that their proofs require the use of large cardinals. 
Specifically, we derive a  new family of ZFC independent theorems  closely related to those of~\cite{hf:alc} (\ref{thm:jfhcaplog}).
Next, we derive a family of theorems (\ref{thm:subsumpoly}) structurally almost identical to and following easily from the theorems of~\ref{thm:jfhcaplog}. 
We make the natural conjecture that the theorems  of the family~\ref{thm:subsumpoly} cannot be proved in ZFC.
We show, however, that a large subclass of the theorems~\ref{thm:subsumpoly} 
  follow from  the statement "subset sum is solvable in polynomial time."  Thus, if our  conjecture that the theorems of~\ref{thm:subsumpoly} can't be proved in ZFC is true, "subset sum is solvable in polynomial time" cannot be proved in ZFC.  We interpret this curious connection between the theory of large cardinals and the $P$ vs $NP$ problem as indicating how difficult it would be to give a ZFC proof of $P = NP$.
\end{abstract}

\maketitle
\section{Introduction}
Basic references are Friedman's 1998 Annals of Mathematics paper, 
{\em Finite functions and the necessary use of large cardinals}
~\cite{hf:nlc}, and his downloadable manuscript
{\em Applications of large cardinals to graph theory}~\cite{hf:alc}.
In Sections 2, 3 and 4 we develop background material and intuition related to certain recursively constructed families of functions on finite subsets of $N^k$, $N$ the 
nonnegative integers. In Section 5, 
we extend  a technique of Friedman\cite{hf:alc},  Theorem 3.4 in addition to Theorem~4.4 and Theorem~4.15,  for creating new independent combinatorial results related to his ZFC independent Jump Free Theorem.  In Section 6, we use these results to relate the classical subset sum problem to the techniques developed in Section 5. 

\section{Elementary background}
We denote by $N$ the set of all nonnegative integers. 
For $z=(n_1, \ldots, n_k)\in N^k$, $\max\{n_i\mid i=1,\ldots, k\}$ is
denoted by $\max(z)$.  We define $\min(z)$ similarly. 
\begin{defn}[\bf Downward directed graph]
\label{def:dwndrctgrph}
We denote by $G=(N^k,\Theta)$ a directed graph (vertex set $N^k$, edge set 
$\Theta$).
If every $(x,y)$ of $\Theta$ satisfies $\max(x) >\max(y)$  then we call $G$ a {\em downward directed lattice graph}. 
For $z\in N^k$, let $G^z = \{x: (z,x)\in \Theta\}$ denote the vertices of $G$ {\em adjacent} to $z$. All such $G=(N^k,\Theta)$ that we consider will be {\em downward directed.}\\
\end{defn}
\begin{defn} [\bf Vertex induced subgraph $G_D$]
For $D\subset N^k$  let 
$G_D = (D, \Theta_D)$ be the subgraph of $G$ with vertex set $D$ and edge set 
$\Theta_D =\{(x,y)\mid (x,y)\in \Theta,\, x, y \in D\}$. We call $G_D$ the {\em subgraph of $G$ induced by $D$}. 
\label{def:vertexinduced}
\end{defn}
\begin{defn}[\bf Cubes and Cartesian powers in $N^k$]
The set  $E_1\times\cdots\times E_k$, where $E_i\subset N$, $|E_i|=p$, $i=1,\ldots, k,$   is called a $k$-cube of length $p$.  If $E_i = E, i=1,\ldots, k,$ then this cube is  $E^k$, the $k$th Cartesian power of $E$.
\label{def:cubespowers}
\end{defn}
\begin{defn}[\bfseries Equivalent ordered $k$-tuples]
\label{def:ordtypeqv}
Two k-tuples in $N^k$, $x=(n_1,\ldots,n_k)$ and $y=(m_1,\ldots,m_k)$, are  
{\em order equivalent tuples $(x\, ot \,y)$} if 
$\{(i,j)\mid n_i < n_j\} =  \{(i,j)\mid m_i < m_j\}$ and  $\{(i,j)\mid n_i = n_j\} =  \{(i,j)\mid m_i =m_j\}.$  
\end{defn}
Note that $ot$ is  an equivalence relation on $N^k$.
The standard SDR (system of distinct representatives) for $ot$  
consists of $({r}_{S_x}(n_1),\ldots,\mathbf{r}_{S_x}(n_k))$ where
 $\mathbf{r}_{S_x}(n_j)$ is the {\em rank} of $n_j$ in $S_x = \{n_1, \ldots, n_k\}$ (e.g,
$x=(3, 8, 5, 3, 8)$, $S_x = \{x\} = \{3,8,5,3,8\}=\{3, 5, 8\}$, $\mathbf{r}(x)=(0, 2, 1, 0, 2)$).
The number of equivalence classes is $\sum_{j=1}^k \sigma(k, j)< k^k,$ $k\geq 2$,
where $\sigma(k,j)$ is the number of surjections from a $k$ set to  a $j$ set.
We use ``$x\,ot\,y$'' and ``$x,\,y$ of order type $ot$'' to mean $x$ and $y$ belong to the same order type equivalence class.

\section{Basic definitions and theorems}

We present some basic definitions due to Friedman~\cite{hf:alc}, \cite{hf:nlc}.

\begin{defn}[\bf Regressive value]
\label{def:regval}
Let $Y\subseteq N$, $X\subseteq N^k$ and $f:X\rightarrow Y$.  An integer $n$
is a  {\em regressive value} of $f$ on $X$
if there exist $x$ such that $f(x)=n<\min(x)$ .
\end{defn}

\begin{defn}[\bf Field of a function and reflexive functions]
\label{def:fldreflfncs}
For $A\subseteq N^k$ define ${\rm field}(A)$ to be the set of all coordinates of elements of $A$.  A function $f$ is {\em reflexive} in $N^k$ if 
${\rm domain}(f) \subseteq N^k$ and  ${\rm range}(f) \subseteq {\rm field}({\rm domain}(f))$.
\end{defn}

\begin{defn}[{\bf The set of functions} $T(k)$ ]
\label{def:tk}
$T(k)$ denotes all reflexive functions with finite domain.
\end{defn}

\begin{figure}[h]
\begin{center}
\includegraphics[scale=.85]{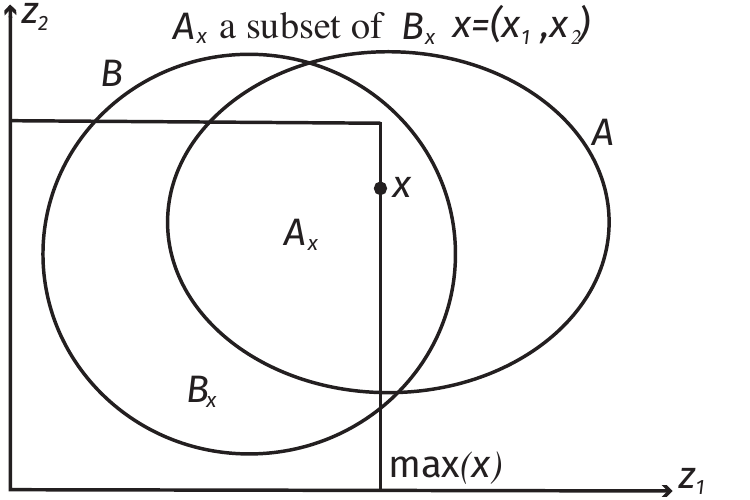}
\caption{Basic jump free condition~\ref{def:fullrefjf}.}
\label{fig:jfvenn}
\end{center}
\end{figure}

\begin{defn} [\bf Full and jump free families] 
Let $Q\subseteq T(k).$
\begin{enumerate}

\item {\bf full family:}  We say that $Q$ is a {\em full} family of functions on $N^k$ if for every finite subset 
$D\subset N^k$ there is at least one function $f$ in $Q$ whose domain is $D$.

\item{\bf jump free family:} For $D\subset N^k$ and $x\in D$ define $D_x = \{z\mid z\in D,\, \max(z) < \max(x)\}$. 
Suppose that for all $f_A$ and $f_B$  in $Q$, where $f_A$ has domain $A$ and $f_B$ has domain $B$,  the conditions
 $x\in A\cap B$, $A_x \subseteq B_x$, and $f_A(y) = f_B(y)$ for all $y\in A_x$ imply that 
$f_A(x) \geq f_B(x)$.  Then $Q$ will be called a {\em jump free} family of functions on $N^k$ (see Figure \ref{fig:jfvenn}). 
\end{enumerate}
\label{def:fullrefjf}
\end{defn} 

\begin{defn}[\bf Function regressively regular over $E$]
\label{def:regreg}
Let $k\geq 2$, $D\subset N^k$, $D$ finite, $f: D\rightarrow N$. 
We say that $f$ is {\em regressively regular} over 
$E$, $E^k\subset D$, if for each order type equivalence class $ot\,$ of $k$-tuples of $E^k$ either (1) or (2) occurs:
\begin{enumerate}
\item{\bf constant less than min $\mathbf{E}$:}  For all $x, y\,\in E^k$ of order type $ot$, $f(x)=f(y)< \min(E).$ 
\item{\bf greater or equal min:}   For all $x\in E^k$ of order type $ot$ $f(x)\geq \min(x).$
\end{enumerate}
\end{defn}   

We use ZFC for the axioms of set theory: Zermelo-Frankel  plus the axiom of choice. 
The Jump Free Theorem, Theorem~\ref{thm:jumpfree} stated below, can be proved in 
 ZFC + ($\forall n$)($\exists$ $n$-subtle cardinal)
 but not in ZFC+ ($\exists$ $n$-subtle cardinal) for any fixed $n$ (assuming this theory is consistent).
A proof of the Jump Free Theorem is in Section 2 of \cite{hf:alc},
``Applications of Large Cardinals to Graph Theory.''  This proof references certain results from~\cite{hf:nlc}.

\begin{thm}[\bf Jump Free Theorem (\cite{hf:alc}, \cite{hf:nlc})] 
\label{thm:jumpfree}
Let $p, k\geq 2$ and $S\subseteq T(k)$ be  a full and jump free family.
Then some $f\in S$ has at most $k^k$ regressive values on some 
$E^k \subseteq {\rm domain}(f)$, $|E| = p$.  
In fact, some $f\in S$ is regressively regular over some $E$ of cardinality $p$.
\end{thm}
We note that the statement ``$f\in S$ is regressively regular over some $E$ of cardinality $p$'' implies that  ''some $f\in S$ has at most $k^k$ regressive values on some $E^k \subseteq {\rm domain}(f)$, $|E| = p$''.  We sometimes keep both statements to emphasize the important $k^k$ bound.

Intuitively, referring to Figure~\ref{fig:jfvenn}, suppose that the region $A_x$ is to be searched for the smallest of some quantity and the result recorded at $x$.  Next, the search region is expanded to a superset $B_x$ with the search results for $A_x$ still valid 
(i.e., $f_A(y) = f_B(y)$ for all $y\in A_x$).  Then, clearly $f_A(x) \geq f_B(x)$.   This expansion property of search algorithms occurs, perhaps somewhat disguised, in many examples.


We next discuss a class of geometrically natural problems that give rise to applications of the Jump Free Theorem. 
Using standard terminology, we use $(x_1, \ldots, x_s)$ to denote a directed path of length $s$ in $G_D$.  If $z\in D$, $(z)$ denotes a path of length one.  A path
$(x_1, \ldots, x_s)$ is terminal if $G_D^{x_s} = \emptyset$.
\begin{figure}[h]
\begin{center}|
\includegraphics[scale=.75]{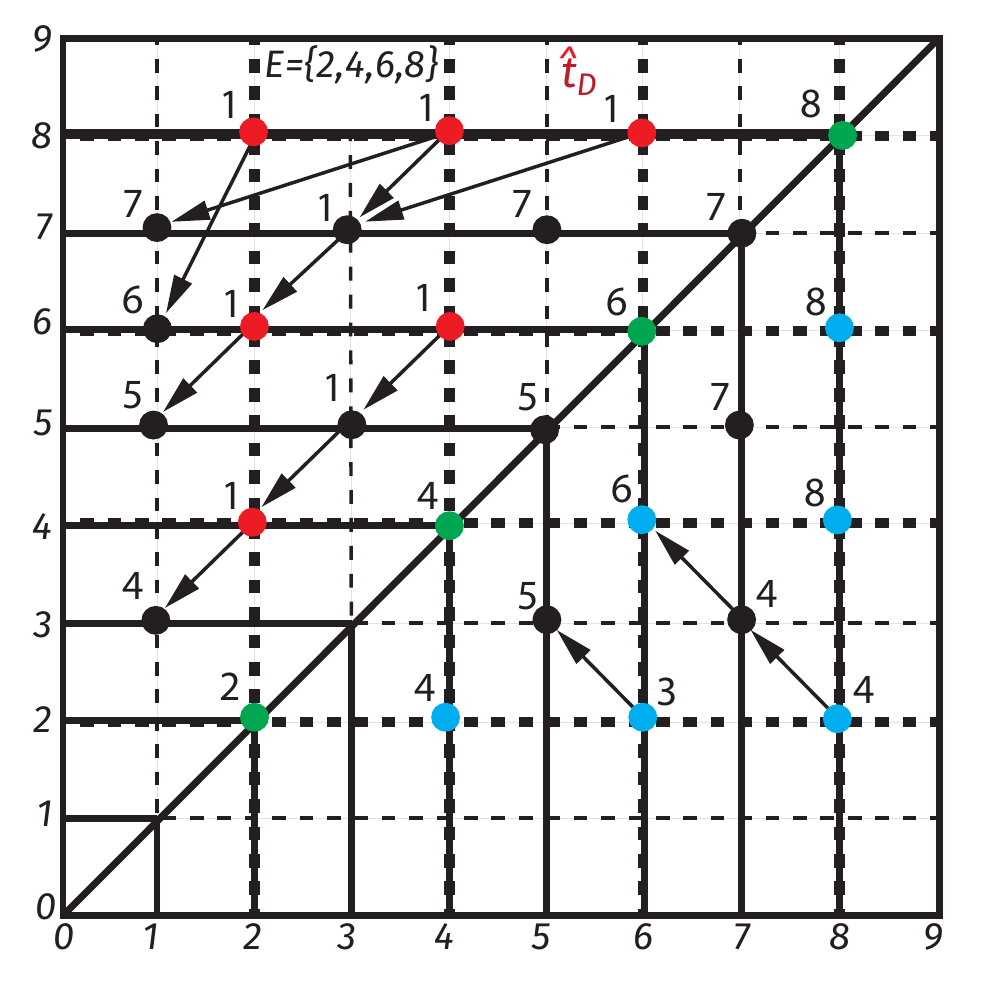}
\caption{$\hat{t}_D$ regressively regular over $E=\{2,4,6,8\}$}
\label{fig:tDhat}
\end{center}
\end{figure}

\begin{defn} [\bf $\hat{t}_D$ terminal path label function]
\label{def:termlabel}
For finite $D\subset N^k$, let $G_D = (D, \Theta_D)$.  Let $T_D(z)$ be the set of all last vertices of terminal paths $(x_1, x_2, \ldots, x_t)$
where $z=x_1$.
Define $\hat{t}_D$ (domain $D$, range ${\rm field}(D)$) by

\begin{enumerate}
\item $\hat{t}_D(z) = \max(z)\;\mathrm{if}\;(z)\;\mathrm{terminal, otherwise}$ 
\item $\hat{t}_D(z)  = \min(\{\min(x)\mid x\in T_D(z)\}).$
\end{enumerate}

We call $\hat{t}_D$ the 
{\em terminal path label function}.
\end{defn}

The choice $\hat{t}_D(z) = \max(z)$ for terminal $(z)$
is used instead of the more natural $\hat{t}_D(z) = \min(z).$ 
This convention makes possible the following application of the Jump Free Theorem
(due to Friedman~\cite{hf:alc}).

\begin{lem}[\bfseries $\{\hat{t}_D\}$ full, reflexive, jump free]
\label{lem:hattfrjf}
Consider the family of functions
$S=\{\hat{t}_D\mid D\subset N^k, D\,\,\mathrm{finite}\}.$ 
Then $\,S$ is full, reflexive, and jump free.

\begin{proof} 
Full and reflexive is immediate.
By the downward condition, $\hat{t}_D=\max(z)$ if and only if $(z)$ is terminal (i.e., $G^z_D = \emptyset$).
Let $\hat{t}_A$ and $\hat{t}_B$  satisfy the conditions for $f_A$ and $f_B$
in Definition~\ref{def:fullrefjf} (2). 
Notice that by definition, $x\notin A_x \,{\rm and}\; x\notin B_x$.
If $(x)$ is terminal in $A$ then $\hat{t}_A(x) =\max(x)\geq \hat{t}_B(x)$ by the downward condition on $G$.
Else, let $(x, \ldots, y)$ be a terminal path in $G_A$.  
Then $\hat{t}_B(y)=\hat{t}_A(y)=\max(y)$ implies  $(x, \ldots, y)$ is a
terminal path in $G_B$. Thus, $\hat{t}_A(x) \geq \hat{t}_B(x)$ as desired.
\end{proof}
\end{lem}

\begin{thm}[\bfseries Jump Free Theorem for $\hat{t}_D$]
\label{thm:jfthat}
Let
$S=\{\hat{t}_D\mid D\subset N^k, D\,\,\mathrm{finite}\}$ and
let $p, k\geq 2$.   Then some $f\in S$ has at most $k^k$ regressive values on some $E^k \subseteq {\rm domain}(f)$, $|E|=p$.  In fact, some $f\in S$ is regressively  regular over some $E$ of cardinality $p$.

\begin{proof}
Follows from Lemma~\ref{lem:hattfrjf} and the Jump Free Theorem~\ref{thm:jumpfree}.  
\end{proof}
\end{thm}

Figure~\ref{fig:tDhat} shows an example of $\hat{t}_D$ regressively regular over a set $E=\{2,4,6,8\}$, where $D\subset N^2$, $|D| = 28$.
Theorem~\ref{thm:jfthat} is one of the most simple combinatorial results in what we call {\em ZFC limbo} -- has a proof using large cardinal assumptions but no known proof within ZFC  itself.

We discuss more complex generalizations in the next section.

\section{More general recursive constructions}

\begin{defn} [{\bf Partial selection}]
\label{def:partselect}
A function $F$ with domain a subset of $X$ and range a subset of $Y$ will be called a {\em partial function}
from $X$ to $Y$ (denoted by $F: X\rightarrow Y$).  If $z\in X$ but $z$ is not in the domain of $F$, we say 
$F$ is {\em not defined at} $z$.
Let $r
\geq 1$.  A partial function 
$F: N^k\times(N^k \times N)^r \rightarrow N$
will be called a {\em partial selection} function ~\cite{hf:alc} if whenever 
$F[x, (y_1,n_1), (y_2,n_2), \ldots (y_r,n_r)]$ is defined we have 
$F[x, (y_1,n_1), (y_2,n_2), \ldots (y_r,n_r)] = n_i$ for some $1\leq i \leq r$.
\end{defn}

\begin{defn}[\bf max constant sets $D_a$]
\label{def:maxcon}
Let $N^k\supset D,\,$ $D$ finite. Let $D_a=\{x\mid x\in D,\, \max(x)=a\}$.
Let $m_0 < m_1 < \cdots < m_q$ be the integers $n$ such that $D_n\neq \emptyset$.
\end{defn}

\begin{defn} [\bf Committee model {\bf $\hat{s}_D$} \cite{hf:alc, gw:lem}]
\label{def:chanlabel}
Let $r\geq 1,$ $k\geq 2$, $G=(N^k,\Theta),$ $G_D = (D, \Theta_D)$, $D$ finite, $G_D^z= \{x\mid (z,x)\in \Theta_D\}.$  
Let $F: N^k\times(N^k \times N)^r \rightarrow N$ be a partial selection function.
If $G_D^z = \emptyset$ define $\Phi^D_z=\emptyset.$
Thus, $\Phi^D_z=\emptyset$ if $z\in D_{m_0}$.
 We define $\Phi^D_z$ and $\hat{s}_D(z)$  (domain $D$, range ${\rm field}(D)$) recursively (on the $m_t$, $t=0, \dots, q$)  as follows. Let
\[
\Phi^D_z = \{ F[z, (y_1,n_1), (y_2,n_2), \ldots, (y_r,n_r)],\;y_i \in G^z_D\}
\]
be the set of defined values of $F$  where  
$n_i=\hat{s}_D(y_i)$ if $\Phi^D_{y_i}\neq\emptyset$ and
$\;n_i=\min(y_i)$ if $\Phi^D_{y_i}=\emptyset.\;$
If $\Phi^D_z\neq\emptyset$,  define $\hat{s}_D(z)$ to be the minimum over $\Phi^D_z$.  
If $\Phi^D_z=\emptyset$,  we set $\hat{s}_D(z) = \max(z)$.
\end{defn}

{\bf NOTE:}   An easy induction on $\max(z)$ shows 
 $\hat{s}_D(z) \leq \max(z)$ with equality if and only if $\Phi^D_z=\emptyset$. 
We give a proof and introduce some terminology.

\begin{lem}[\bf $\hat{s}_D(z)$ structure ]
 $\hat{s}_D(z)\leq \max(z)$ with
 $\hat{s}_D(z) = \max(z)$ if and only if $\Phi^D_z=\emptyset$.
 \begin{proof}
 We use induction on $\max(z)$ to construct $\hat{s}_D(z)$ and  $\Phi^D_z$.
Let $D_a=\{x\mid x\in D,\, \max(x)=a\}$.
Let $m_0 < m_1 < \cdots < m_q$ be the list of $n$ such that $D_n\neq \emptyset$.
If $z\in D_{m_0}$ then the set of adjacent vertices $G_D^z = \emptyset$.
Thus, $\Phi^D_z = \emptyset$ and $\hat{s}_D(z) = \max(z)$ for all $z\in D_{m_0}$.
In general, assume that for $t<j$, $z\in D_{m_t}$, $\hat{s}_D(z)\leq \max(z)$ with
 $\hat{s}_D(z) = \max(z)$ if and only if $\Phi^D_z=\emptyset$.
 Consider $z\in D_{m_j}$.
If (1)  $\Phi^D_z = \emptyset$ then $\hat{s}_D(z) = \max(z).$
If (2)  $\Phi^D_z \neq \emptyset$  let $n=F[z, (y_1,n_1), (y_2,n_2), \ldots, (y_r,n_r)]\in \Phi^D_z $, $\;y_i \in G^z_{D}$  thus $y_i\in D_t$, $t<j$.

First, if $\Phi^D_{y_i}=\emptyset.\;$ then  $\;n_i=\min(y_i)< \max(z)$. 

Second,  if $\Phi^D_{y_i}\neq\emptyset$ 
then, by the induction hypothesis, $n_i=\hat{s}_D(y_i)<\max(y_i)<\max(z)$.   
 Thus, $\hat{s}_D(z)\leq \max(z)$ with
 $\hat{s}_D(z) = \max(z)$ if and only if $\Phi^D_z=\emptyset$.
 \end{proof}
 \end{lem}
 
The following result is due to Friedman~\cite{hf:alc}.

\begin{thm}[\bfseries Large scale regularities for $\hat {s}_D$ (\cite{hf:alc}]
\label{thm:jfhats}
Let $r\geq 1$, $p, k\geq 2$.
$S=\{\hat {s}_D\mid D\subset N^k, D\,\mathrm{finite}\}$. Then some 
$f\in S$ has at most $k^k$ regressive values on some 
$E^k \subseteq {\rm domain}(f)$, $|E|=p.$ 
In fact, some $f\in S$ is regressively regular over some $E$ of cardinality $p$.

\begin{proof}
Recall the Jump Free Theorem~\ref{thm:jumpfree}. 
Let  $S=\{\hat {s}_D\mid D\subset N^k,\;D\,\mathrm{finite}\}.$
$S$ is obviously full and reflexive. 
We show that $S$ is jump free.
We show for all $\hat {s}_A$ and $\hat {s}_B$  in $S$,  the conditions
 $x\in A\cap B$, $A_x \subseteq B_x$, and $\hat {s}_A(y) =\hat {s}_B(y)$ for all $y\in A_x$ imply that 
$\hat {s}_A(x) \geq \hat {s}_B(x)$  (i.e., $S$ is  {\em jump free}). 
If $\Phi^A_x = \emptyset$ then $\hat {s}_A(x)=\max(x)\geq \hat {s}_B(x).$

Assume $\Phi^A_x \neq \emptyset.$ 
Let $n=F[x, (y_1,n_1), (y_2,n_2), \ldots (y_r,n_r)]\in \Phi^A_x$ 
(observe that $y_i\in G_A^x \subseteq G_B^x$) where 
$n_i=\hat {s}_A(y_i)$ if $\hat {s}_A(y_i)<\max(y_i)$ 
(i.e., $\Phi^A_{y_i} \neq \emptyset$) and 
$n_i=\min(y_i)$ if $\hat {s}_A(y_i)=\max(y_i)$ (i.e., $\Phi^A_{y_i} = \emptyset$).
But $\hat {s}_A(y_i) =\hat {s}_B(y_i)$, $i=1, \ldots, r,$ implies
$n\in \Phi^B_x$ and thus $\Phi^A_x \subseteq \Phi^B_x$ and
$
\hat {s}_A(x)=\min(\Phi^A_x) \geq \min(\Phi^B_x)=\hat {s}_B(x).
$
\end{proof}
\end{thm}

Next we give  an example of $\hat{s}_D$.
 
 \begin{figure}[h]
\begin{center}
\includegraphics[scale=.7]{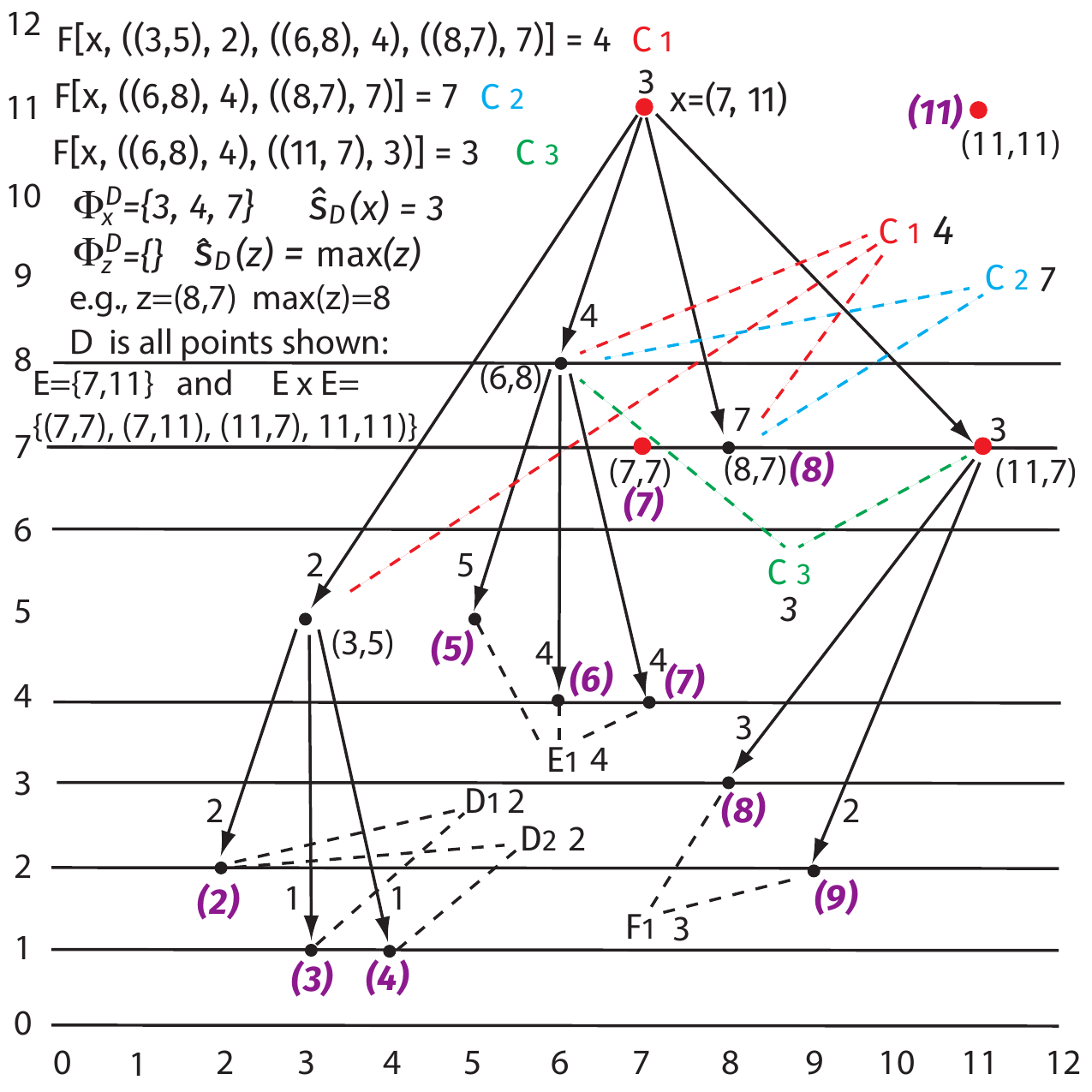}
\caption{An example of $\hat{s}_D$}
\label{fig:univs}
\end{center}
\end{figure}

As an example of computing $\hat{s}_D$, consider Figure~\ref{fig:univs}.
The computation is recursive on the $\max$ norm (and doesn't illustrate all of the subtleties).  
The $\hat{s}_D$ values of the non-isolated terminal vertices where  $\Phi^A_x = \emptyset$ are shown in parentheses,  left to right: (2), (3), (4), (5), (6), (7), (8), (8), (9).  These numbers are  $\max((a,b))$ for each such terminal vertex $(a,b)$.
Partial selection functions are of the form 
$F:N^2 \times (N^2\times N)^r \rightarrow N$ ($r=2, 3$ here).  
In particular we have
$F[x, ((3,5),2), ((6,8),4), ((8,7),7)]=4,\;\;$
$F[x, ((6,8),4), ((8,7),7)]= 7$, and
$F[x, ((6,8),4), ((11,7),3)]= 3$.
Intuitively, we think of these as (ordered) committees reporting values to the boss, $x=(7,11).$
The first committee, $\rm{C}1$, consists of subordinates, $(3,5),(6,8), (8,7)$
reporting respectively $2, 4, 7$.  
The committee decides to report $4$ (indicated by  $\rm{C}1\;4$ in  
Figure~\ref{fig:univs}). 
The recursive construction starts with terminal vertices reporting their minimal coordinates.
But, the value reported by each committee is not, in general, the actual minimum of the reports of the individual members. 
Nevertheless, the boss, $x=(7,11)$,  {\em always} takes the minimum of the values reported by the committees.  
In this case the values reported by the committees are $4, 7, 3$ the boss takes $3$ (i.e., $\hat{s}_D(x)=3$ for the boss, $x=(7,11)$).
Observe that a function like $F((7,11), ((6,8),4), ((8,7),7)$ where $r=2$, can be padded
to the case $r=3$
(e.g., $F((7,11), ((6,8),4), ((8,7),7), ((8,7),7))).$ 

Observe in Figure~\ref{fig:univs} that the values in parentheses, (2), (3), (4), (5), (6), (7), (8), (8), (9), don't figure into the recursive construction of $\hat{s}_D$.  They immediately pass their minimum values on to the computation: 2, 1, 1, 5, 4, 4, 7, 3, 2.
Here, $E^2 = \{(7,7), (7,11), (11,7), (11,11)\}$. Isolated terminal vertices are
$\{(7,7),(11,11)\}$, the diagonal of $E^2$.

\section{Combinatorial Generalizations}
In this section we present some results that are based on results of Friedman~\cite{hf:alc} (specifically, Theorem~4.4 and the ideas of Theorem~4.1 and the earlier Theorem~3.3).  Friedman removes any mention of the graph $G$ and works with an ``equivalent streamlined version.''  By contrast we focus on the graph model in this discussion.

We extend Friedman's results slightly by introducing a class of functions 
$\mathbb{R}=\{\rho_{D}\mid D\subset N^k,\; D\;{\rm finite}\,, \rho_D:D \rightarrow N,\, \min(x)\leq \rho_{D}(x),\,x\in D\}$. These ``min dominant'' functions allow us to relax the reflexive condition. 
We follow closely Definition~\ref{def:chanlabel}.
\begin{defn}[{\bf $h^{\rho_D}$ initialized by $\mathbb{R}$}]
\label{def:genhat}
Let $r\geq 1,$ $k\geq 2$, $G=(N^k,\Theta),$ $G_D = (D, \Theta_D)$, $D$ finite, $G_D^z= \{x\mid (z,x)\in \Theta_D\}.$ 
Let $F: N^k\times(N^k \times N)^r \rightarrow N$ be a partial selection function.
An {\em initializing min dominant family of functions} is specified as follows:
\[
\mathbb{R}=\{\rho_{D}\mid N^k\supset D\; {\rm finite}\,, \rho_D:D \rightarrow N,\, \rho_{D}(x)\geq \min(x),\,x\in D\}.
\]
We define $\Phi^{\rho_D}_z$, $\,h^{\rho_D}$ recursively on $\max(z)$ (using the $m_t$, $t=0, \dots q$ of Definition~\ref{def:maxcon}) . 
If $G_D^z = \emptyset$ define $\Phi^{\rho_D}_z = \emptyset$. Thus, $\Phi^{\rho_D}_z = \emptyset,\;z\in D_{m_0}.\;$
We define $\Phi^{\rho_D}_z$ and $h^{\rho_D}(z)$ recursively   as follows.
Let 
\[
\Phi^{\rho_D}_z = \{ F[z, (y_1,n_1), (y_2,n_2), \ldots, (y_r,n_r)],\;y_i \in G^z_D\}
\]
be the set of defined values of $F$  where  
$n_i=h^{\rho_D}(y_i)$ if $\Phi^{\rho_D}_{y_i}\neq\emptyset$ and
$\;n_i=\min(y_i)$ if $\Phi^{\rho_D}_{y_i}=\emptyset.\;$
If $\Phi^{\rho_D}_z=\emptyset$,  define $h^{\rho_D}(z) = \rho_{D}(z)$.
If $\Phi^{\rho_D}_z\neq\emptyset$,  define $h^{\rho_D}(z)$ to be the minimum over 
$\Phi^{\rho_D}_z$.
We say {\bf $h^{\rho_D}$ is initialized by $\mathbb{R}$}.
\end{defn}

Recall Definition~\ref{def:chanlabel} and the recursive construction of
 $\hat{s}_D$ and $h^{\rho_D}$.  
 Note that $h^{\rho_D}$ is the same as $\hat{s}_D$ if $\rho_D(x) = \max(x)$.
 The inductive structure of the following lemma is ``intuitively obvious'' but we give a formal proof anyway.
\begin{lem}[{\bf Compare} $\hat{s}_D,\,h^{\rho_D}$] 
\label{lem:shatvsh}
For all $z\in D$, either (1) $\Phi^D_z=\Phi^{\rho_D}_z = \emptyset$ 
 and
$h^{\rho_D}(z) = \rho_D(z)$, $\;\hat{s}_D(z)=\max(z)\,$ or (2)
$\,\Phi^D_z=\Phi^{\rho_D}_z\neq \emptyset$ and $h^{\rho_D}(z)=\hat{s}_D(z)$.
\begin{proof}
We define $D_a=\{x\mid x\in D,\, \max(x)=a\}$.
Let $m_0 < m_1 < \cdots < m_q$ be the integers $n$ such that $D_n\neq \emptyset$.
If $z\in D_{m_0}$ then $G_D^z = \emptyset$.
Thus, $\Phi^D_z=\Phi^{\rho_D}_z = \emptyset$  and $h^{\rho_D}(z) = \rho_D(z)$, 
$\hat{s}_D(z)=\max(z),\,z\in D_{m_0}$.

Assume, for all $z\in D_{m_t}$, $0\leq t<j$, either $\Phi^D_z=\Phi^{\rho_D}_z = \emptyset$  and $h^{\rho_D}(z) = \rho_D(z)$, $\;\hat{s}_D(z)=\max(z)$
or $\Phi^D_z=\Phi^{\rho_D}_z\neq \emptyset$ and 
$h^{\rho_D}(z)=\hat{s}_D(z)$. 

Let $z\in D_{m_{j}}$.
If $\Phi^D_z=\Phi^{\rho_D}_z = \emptyset$, then $h^{\rho_D}(z) = \rho_D(z)$ and 
$\hat{s}_D(z)=\max(z)$.
Otherwise, either $\Phi^D_z\neq \emptyset$ or $\Phi^{\rho_D}_z\neq\emptyset$.
Assume that $\Phi^{\rho_D}_z\neq\emptyset$. 
Let 
\[
\Phi^{\rho_D}_z = \{ F[z, (y_1,n_1), (y_2,n_2), \ldots, (y_r,n_r)],\;y_i \in G^z_D\}.
\]
Choose $n=F[z, (y_1,n_1), (y_2,n_2), \ldots, (y_r,n_r)]\in \Phi^{\rho_D}_z$. 
Thus  $y_i\in D_{m_t}$ for some $t<j$.
By the induction hypothesis, either 
$(1)\;\Phi^D_{y_i}=\Phi^{\rho_D}_{y_i} = \emptyset$, $\;h^{\rho_D}({y_i}) = \rho_D({y_i})\,$ and $\,\hat{s}_D({y_i})=\max({y_i})$, in which case $n_i = \min({y_i})$, or
$(2)\;\Phi^D_{y_i}=\Phi^{\rho_D}_{y_i}\neq \emptyset$ and $h^{\rho_D}({y_i})=\hat{s}_D({y_i})=n_i$.   In either case, $n\in \Phi^D_z$ and thus 
$\Phi^{\rho_D}_z \subseteq \Phi^D_z$. 
In the same manner we conclude that $ \Phi^D_z \subseteq \Phi^{\rho_D}_z$.
Thus, in fact, $\,\Phi^D_z=\Phi^{\rho_D}_z\neq \emptyset$ and $h^{\rho_D}(z)=\hat{s}_D(z)$.
\end{proof}
\end{lem}
Next, we consider {\em regressive regularity}. 

\begin{lem}[\bf Compare regressive regularity of  $\hat{s}_D,\,h^{\rho_D}$ ]
\label{lem:comregreg}
Let $E$ be of cardinality $\,p\geq 2$.
Then $\hat{s}_D$ is regressively  regular over $E$ if and only if $\,h^{\rho_D}$ is regressively  regular over $E$. 

\begin{proof}

For $z\in D$ we have shown (Lemma~\ref{lem:shatvsh}) there are two cases:
$$(1)\;\;\Phi^D_{z}=\Phi^{\rho_D}_{z} = \emptyset,\; h^{\rho_D}({z}) = \rho_D({z})\; {\rm and} \;\hat{s}_D({z})=\max({z})$$ 
and
$$(2)\;\;\Phi^D_{z}=\Phi^{\rho_D}_{z}\neq \emptyset\; {\rm and} \;h^{\rho_D}({z})=\hat{s}_D({z}).$$ 
 
{\bf First}, we show for all $x, y\,\in E^k$ of order type $ot$, 
$\hat{s}_D(x)=\hat{s}_D(y)< \min(E)$ 
 if and only if $h^{\rho_D}(x)=h^{\rho_D}(y)< \min(E)$. Case 
 (1) above is ruled out because $h^{\rho_D}({z}) = \rho_D({z})\geq \min(z)\geq\min(E)$
and $\hat{s}_D({z})=\max({z})\geq\min(E)$.
Thus we have case (2) $\,\Phi^D_{z}=\Phi^{\rho_D}_{z}\neq \emptyset$ and $h^{\rho_D}({z})=\hat{s}_D({z})$ for $z=x,\, y$.
Thus,  $\hat{s}_D(x)=\hat{s}_D(y)< \min(E)$ if and only if $h^{\rho_D}(x)=h^{\rho_D}(y)< \min(E)$.

{\bf Second}, suppose that for all $x\in E^k$ of order type $ot$, 
$h^{\rho_D}(x)\geq \min(x).$
This set of order type $ot$ splits naturally into two sets, $\{x \mid \Phi^D_x\neq 
\emptyset\}$ and
$\{x \mid \Phi^D_x = \emptyset\}$. 
On the first set, $\min(x) \leq h^{\rho_D}(x) =  \hat{s}_D(x)$ and on the second set $h^{\rho_D}(x) = \rho_D(x)\geq \min(x)$ and $\hat{s}_D(x) = \max(x) \geq \min(x)$.
Thus, $\hat{s}_D(x)\geq \min(x)$.
The same argument works if we assume for $x\in E^k$ of order type $ot$ 
$\hat{s}_D(x)\geq \min(x).$
Thus, for $x\in E^k$ of order type $ot$, $h^{\rho_D}(x)\geq \min(x)$ if and only if
$\hat{s}_D(x)\geq \min(x).$
\end{proof}
\end{lem}

\begin{thm}[\bfseries Regressive regularity of $h^{\rho_D}$]
\label{thm:jfh} 
Let $G=(N^k, \Theta)$, $r\geq 1$, $p, k\geq 2$. 
Let $S=\{h^{\rho_D}\mid D\subset N^k,\; D\; \mathrm{ finite}\}$.  Then some $f\in S$ has at most $k^k$ regressive values on some $E^k \subseteq {\rm domain}(f)=D$, 
$|E|=p$.
In fact, some $f\in S$ is regressively  regular over some $E$ of cardinality~$p$.

\begin{proof}
Follows from Theorem~\ref{thm:jfhats} and Lemmas~\ref{lem:shatvsh}, \ref{lem:comregreg}.
We claim that the set  $S=\{h^{\rho_D}\mid D\subset N^k,\; D\; \mathrm{ finite}\}$ is a full family of functions such that for every $p\geq 2$ there is a function 
$h^{\rho_D}$ which is regressively regular over some $E$, $|E|=p$. Lemmas~\ref{lem:shatvsh} and \ref{lem:comregreg} show that to find such an $E$ for $h^{\rho_D}$ we can invoke
Theorem~\ref{thm:jfhats} and find such an $E$ for $\hat{s}_D$.
\end{proof}
\end{thm}

{\bf Remark: Independence of the families of Theorem~\ref{thm:jfh}.}
\label{rem:ijfh}
From Theorem~\ref{thm:jfh} the regressive regularity of the families of functions
$\{h^{\rho_D}\mid D\subset N^k,\;|D|<\infty\}$ is in ZFC limbo as the only proof we have at this point uses large cardinal assumptions and these assumptions cannot be proved in ZFC.
However, Friedman\cite{hf:alc}, has removed these families from limbo.
In particular, it has been shown by Friedman\cite{hf:alc},  Theorem~4.4 through Theorem~4.15, that a special case of Theorem~\ref{thm:jfh} ($\rho_D=\min$) requires the same large cardinals to prove as the Jump Free Theorem.  
Thus, Theorem~\ref{thm:jfh} provides a family of statements independent  of ZFC and parameterized by a choice of an initializing min dominant family of functions:
\[
\mathbb{R}=\{\rho_{D}\mid N^k\supset D\; {\rm finite}\,, \rho_D:D \rightarrow N,\, \rho_{D}(x)\geq \min(x),\,x\in D\}.
\]

\section{Using the $\rho_D$ and the subset sum problem}

\begin{defn}[\bf $D$  capped by $E^k\subset D$]
\label{def:cap}
For $k\geq 2$, $E^k \subseteq D\subset N^k$, let 
$\max(D)=\max\{\max(z): z\in D\}$.  
Let $\setmax(D)=\{z\mid z\in D, \max(z)=\max(D)\}$.
If  $\setmax(D) = \setmax(E^k)$, we say that $D$ is {\em capped by} 
$E^k\subseteq D$ with the {\em cap} defined to be $\setmax(E^k)$.
\end{defn}

Note that if $D$ is {\em capped by} $E^k\subseteq D$ then $D$ determines $E^k$ uniquely in the obvious way.  
An example is shown in Figure~\ref{fig:tDhat} and Figure~\ref{fig:univs}.

The following theorem is analogous to Theorem~\ref{thm:jfh}.

\begin{thm}[\bfseries Regressively regular $h^{\rho_D}$, capped version]
\label{thm:jfhcap}
Let $G=(N^k, \Theta)$, $r\geq 1$, $p, k\geq 2$. Let $S=\{h^{\rho_D}\mid D\subset N^k,\;D\; {\rm finite}\,,|<\infty\}$.    
Some $h^{\rho_D} \in S$ is regressively  regular over some $E$, $|E|=p$,
$E^k\subseteq D$, $D$ capped by $E^k.$
\begin{proof}
From Theorem~\ref{thm:jfh} there is an  $h^{\rho_D}\in S$ that is regressively regular over some $E$, $|E|=p$, $E^k\subseteq D.$
Let $E=\{e_0, \ldots, e_{p-1}\}$.  
Let $D_x=\{z\mid z\in D, \max(z)< \max(x)\}$.
Let $\widehat{D}= D_{e_{p-1}} \cup \setmax(E^k)$ so $\widehat{D}$ is capped by $E^k$.
Using the downward condition on $G_D$ and hence $G_{\widehat{D}}$
we have the restriction $h^{\rho_D} | \widehat{D}$ is regressively regular over $E$.
Note that $h^{\rho_D} | \widehat{D}$ may or may not be equal to the function $h^{\rho_{\widehat{D}}}\in S.$  
But Lemma~\ref{lem:shatvsh} and Lemma~\ref{lem:comregreg} apply in either case.
Thus we conclude that the function $h^{\rho_{\widehat{D}}}\in S$ is also regressively regular over $E$.
\end{proof}
\end{thm}



\begin{defn}[\bf $t$-log bounded]
\label{def:lgbn}
Let $p, k\geq 2$, $t\geq 1$. 
The function $\rho_D$ is 
 $t$-log {\em bounded} over $E^k\subset D$ where $E=\{e_0, \ldots, e_{p-1}\},\,$  if the cardinality
\[
|\{ \rho_D(x) - \min(x) : 0< \rho_D(x) - \min(x) < e_0 k^k,\, x\in E^k\}| \leq t\log_2 (p^k).
\]
In this case, we write $\rho_D\in \rm{LOG}(k,E,p,D,t)$. 
The set 
$$
\mathbb{R}=\{\rho_{D}\mid \rho_D:D \rightarrow N,\, \min(x)\leq \rho_{D}(x),\,x\in D\}
$$
is $t$-log bounded if $\rho_D\in \rm{LOG}(k,E,p,D,t)$ when $D$ is capped by $E^k.$
In this case we write $\mathbb{R}_t$ for $\mathbb{R}.$
\end{defn}

{\bf Discussion of Definition~\ref{def:lgbn}}. 
Recalling that $\rho_D(x)\geq \min(x)$ and 
$\rho_D(x)$ can be arbitrarily large,
we can choose the cardinality $|\{x: \rho_D(x) - \min(x)   \geq e_0 k^k\}|$ large enough to make $\rho_D\in \rm{LOG}(k,E,p,D,t)$. 
We can also choose the $\rho_D(x) - \min(x) \geq e_0k^k$  distinct. We make that general assumption in what follows.  Also observe that the function $\rho_D$ has to be defined initially, before the $h^{\rho_D}$.  
This is possible because if $D$ is capped by $E^k$ the $E^k$ is uniquely defined.
Thus, $\rho_D$ can be defined to satisfy~~\ref{def:lgbn}.  
\begin{thm}[\bfseries Regressive regularity $t$-log bounded case]
\label{thm:jfhcaplog}
Let $G=(N^k, \Theta)$, $r\geq 1$, $p, k\geq 2$. Let $S=\{h^{\rho_D}\mid D\subset N^k,\;D\; {\rm finite}\,\}$
where the initializing set $\mathbb{R}_t$
is $\;t$-log bounded.
Then some $h^{\rho_D} \in S$ is regressively  regular over some  $E$,  $|E|=p$, 
$E^k\subseteq D$, $D$ capped by $E^k$ and $\rho_D\in LOG(k,E,p,D,t)$.
\begin{proof}
Follows from Theorem~\ref{thm:jfhcap} which states that  
some $h^{\rho_D} \in S$ is regressively  regular over some such $E$, $|E|=p$, 
$E^k\subseteq D$, $D$ capped by $E^k.$
From Definition~\ref{def:lgbn}, for each such capped pair $D$ and $E^k$,  $\rho_D$ has already been defined so that $\rho_D\in \rm{LOG}(k,E,p,D,t)$.  
\end{proof}
\end{thm}

Theorems~\ref{thm:jfhcaplog} and~\ref{thm:jfhcap} require the same large cardinal assumptions as Theorem~\ref{thm:jfh} (see Remark following~\ref{thm:jfh}) .

\begin{defn}({\bf $h^{\rho_D}$ partitions $E^k\subseteq D$ into three blocks).}
\label{def:blocks}
Given any $h^{\rho_D}$ we have a natural partition of $E^k\subseteq D$ into three sets:  
$$E^k_0 =\{ x\in E^k: h^{\rho_D}(x) < \min(E)\}$$
$$E^k_1=\{x\in E^k: \min(E)\leq  h^{\rho_D}(x) < \min(x)\}$$
$$E^k_2=\{x\in E^k: \min(x)\leq  h^{\rho_D}(x) \}.$$
\end{defn}

In Definition~\ref{def:setsdis} below we associate sets of integers with each of the three blocks of this partition.  Our associated sets are chosen because of their natural and generic relationship to regressive regularity.
Let $Z =\{0, \pm 1, \pm 2, \ldots\}$ be the integers and let 
$\{I_D\mid I_D: D\rightarrow Z, D\subset N^k,\,D\, {\rm finite}\}$,  be a family of functions defined for all finite subsets of $N^k$.
We use the terminology of Theorem~\ref{thm:jfhcaplog}.

\begin{defn}[\bf Sets of instances]
\label{def:setsdis}
Let $S=\{h^{\rho_D}\mid D\subset N^k,\;D\; {\rm finite}\,\}$
with initializing set $\mathbb{R}_t$. Using Definition~\ref{def:blocks} we define
\[
\Delta h^{\rho_D}E_0^k = \{h^{\rho_D}(x) - \min(E) : x\in E_0^k\}\;\;\;\;\;\;
\]
\[
\Delta h^{\rho_D}E_1^k = \{I_D(x) : x\in E_1^k\}\;\;\;\;\;\;
\]
\[
\Delta h^{\rho_D}E_2^k = \{\rho_D(x) - \min(x) : x\in E^k_2 \}.\;\;\;\;\;\;
\]
\end{defn}

The sets introduced in Definition~\ref{def:setsdis} (to be used as sets of instances in the proof of Theorem~\ref{thm:subsumpoly}) are constructed to be sensitive to the case where $h^{\rho_D}$ is regressively regular over $E$. Observe that $|\cup _{i=0}^2\, \Delta h^{\rho_D}E_i^k|\leq p^k$,  $|E|=p$.

We summarize terminology:
\begin{enumerate}
\item  $N$ nonnegative integers, $Z$ integers.
\item  $N^k$ nonnegative integral lattice, $k\geq 2$. 
\item $\{I_D\mid I_D: D\rightarrow Z,  {\rm finite}\,D\subset N^k\}$ family of functions.
\item  $E=\{e_0, \dots, e_{p-1}\} \subset N$, $|E|=p\geq 2$.
\item  $E^k\subseteq D$, $D$ capped by $E^k$ defined by $\setmax(D) = \setmax(E^k)$.
\item  $\mathbb{R} = \{\rho_D\mid \rho_D:D\rightarrow N, {\rm finite\,} D\subset N^k, \rho_D(x)\geq x \}$
initializing family.
\item $\mathbb{R}_t$ a $t$-log bounded initializing family of functions, $t\geq 1$.
\item  $F:N^k\times (N^k\times N)^r\rightarrow N,\;$ $r\geq 1,\;$ partial selection functions.
\item  $G=(N^k,\Theta)$ downward directed graphs on $N^k$.
\item  $G_D = (D, \Theta_D)$ restriction of $G$ to $D$.
\item  $h^{\rho_{D}}$ functions initialized by $\mathbb{R}$ defined recursively on finite $D\subset N^k$.
\end{enumerate}

\begin{thm}[\bf Subset sum connection]
\label{thm:subsumpoly}
For fixed $k, t, F, G, $ consider sets of instances
of the form
$H^{k,t}_{F,G}(E,p,D) =\cup _{i=0}^2\, \Delta h^{\rho_D}E_i^k$
where 
$ h^{\rho_D}$ is initialized by $\mathbb{R}_t$, $E^k\subseteq D$ , $D$ capped by $E^k$, $|E|=p$ and  the $\Delta h^{\rho_D}E_i^k$
are defined in Definition~\ref{def:setsdis}.
For each $p$ there exists $\hat{E}$ and $\hat{D}$, $|\hat{E}|=p$, such that
the subset sum problem for
\[
\{H^{k,t}_{F,G}(\hat{E},p,\hat{D}):  p=2, 3, \ldots\}
\]
is solvable in time $O(p^{kt})$.  
\begin{proof}
From Theorem~\ref{thm:jfhcaplog}, 
for any $p$, we  can choose ${\hat{D}}$ capped by $\hat{E}^k$, $|\hat{E}|=p$, such that 
$h^{\rho_{\widehat{D}}}$ is regressively regular over  $\hat{E}$.
For notational simplicity we write $\hat{E}=\{e_0, \ldots, e_{p-1}\}.$

By regressive regularity,  
the set $E^k_1=\{x\in \hat{E}^k: \min(\hat{E})\leq  h^{\rho_{\widehat{D}}}(x) < \min(x)\}$ is empty, thus
$\Delta h^{\rho_{\widehat{D}}}\hat{E}_1^k  = \emptyset$ (see Definition~\ref{def:setsdis}).
 
For  $\Delta h^{\rho_{\widehat{D}}}\hat{E}_0^k $ we have $h^{\rho_{\widehat{D}}}(x) - e_0 < 0$.
Observe that  $|h^{\rho_{\widehat{D}}}(x) - e_0  |< e_0$ and, by regressive regularity, 
the cardinality $|\Delta h^{\rho_{\widehat{D}}}\hat{E}_0^k|$ is smaller than $k^k$ thus

\begin{equation}
\label{eq:totsum}
\sum_{x\in\Delta h^{\rho_{\widehat{D}}}\hat{E}_0^k } |h^{\rho_{\widehat{D}}}(x) - e_0  |< e_0 k^k.
\end{equation}
From t-log bounded, we have:
\[
|\{\rho_{\widehat{D}}(x) - \min(x) : 0< \rho_{\widehat{D}}(x) - \min(x) < e_0 k^k,\, x\in {\hat{E}}^k\}| \leq t\log_2 (p^k).
\]
The negative terms in the instance $H^{k,t}_{F,G}(\hat{E},p,\hat{D})$ come from
\[
\Delta h^{\rho_{\widehat{D}}}\hat{E}_0^k = \{h^{\rho_{\widehat{D}}}(x) - \min(\hat{E}) : x\in \hat{E}^k,\; h^{\rho_{\widehat{D}}}(x)< \min(\hat{E}) \}.\;\;\;\;\;\;
\]
The cardinality  $|\Delta h^{\rho_{\widehat{D}}}\hat{E}_0^k| < k^k$.
 
The positive terms in the instance come from
\[
\Delta h^{\rho_{\widehat{D}}}\hat{E}_2^k= \{\rho_{\widehat{D}}(x) - \min(x) : x\in \hat{E}^k, \min(x)\leq  h^{\rho_{\widehat{D}}}(x) \}.\;\;\;\;\;\;
\]
If $\,0\in \Delta h^{\rho_{\widehat{D}}}\hat{E}_2^k$ then the solution is trivial as $0$ is the target.
We use equation~\ref{eq:totsum} to rule out having to consider positive values of 
$\rho_D(x) - \min(x)\geq e_0k^k.$
We can check all possible solutions by comparing the sums of less than $2^{k^k}$ subsets of negative terms with less than 
$$2^{t\log_2(p^k)}=2^{tk\log_2(p)}=p^{kt}$$ 
subsets of positive terms. Thus we can check all possible solutions in $O(p^{kt})$ comparisons.
%
%
\end{proof}

\end{thm}
\section{Conclusions}
We have proved Theorem~\ref{thm:subsumpoly} using
Theorem~\ref{thm:jfhcaplog}. Theorem~\ref{thm:jfhcaplog} cannot be proved in  ZFC alone (for each fixed initiating set $\mathbb{R}_t$).
We know of no other proof.  
Thus, Theorem~\ref{thm:subsumpoly} for each fixed $\mathbb{R}_t$, $t\geq 1$, is in ZFC limbo.
If a ZFC proof could be found that the subset sum problem is solvable in polynomial time $O(n^\gamma)$ where $n$ is the length of the instance ($p^k$ for fixed $k$ here), then that result would prove Theorem~\ref{thm:subsumpoly} for $t=\gamma$ and thus remove that case from limbo by showing that it is provable within ZFC. 
We conjecture, however, that Theorem~\ref{thm:subsumpoly} cannot be proved in  ZFC alone (for each fixed initiating set $\mathbb{R}_t$). 
The basis for this conjecture is that the subset sum problem arises from
Theorem~\ref{thm:jfhcaplog} in a very natural, generic way.
Of course, if our conjecture is true, ``subset sum is solvable in polynomial time'' cannot be proved in ZFC (perhaps because it is false).

Specific examples of sets of instances of the form $\Delta h^{\rho_D}E_1^k$ (Definition~\ref{def:setsdis}) are used in earlier versions of this paper (e.g., \cite{gw:sub} and \cite{gw:lim}). 
The set $\Delta h^{\rho_D}E_1^k = \{I_D(x) : x\in E_1^k(x)\}$ used here is more general, making the possibility of a ZFC proof of Theorem~\ref{thm:subsumpoly} intuitively even more difficult.

{\bf Acknowledgments:}  The author thanks  Professor Sam Buss (University of California San Diego, Department of Mathematics), 
Professor Emeritus Rod Canfield (University of Georgia, Department of Computer Science), and Professor Emeritus Victor W. Marek (University of Kentucky, Department of Computer Science) for numerous helpful suggestions.

%

\bibliographystyle{alpha}
\bibliography{jeffjoc}

\end{document}